\documentclass[10pt,twoside]{siamltex}
\usepackage{amsfonts,ifthen,stmaryrd,url}
\usepackage[leqno]{amsmath}

\def\cvd{~\vbox{\hrule\hbox{%
     \vrule height1.3ex\hskip0.8ex\vrule}\hrule } }

\newcommand{\reals}{\mathbb{R}}
\newcommand{\complex}{\mathbb{C}}

\DeclareMathOperator{\tr}{tr}
\DeclareMathOperator{\nz}{nz}
\DeclareMathOperator{\nnz}{nnz}
\DeclareMathOperator{\rem}{rem}

\let\lkq=\obslash
\let\rkq=\oslash

\def\makelenviron#1#2{%
\newenvironment{l#1}[1]{%
 \begin{#1}\ifthenelse{\equal{##1}{}}{}{\label{##1}}#2\ignorespaces%
 }{\end{#1}}%
}

\makelenviron{definition}{\upshape}
\makelenviron{theorem}{}

\newif\ifhavecvd
\let\elacvd=\cvd
\def\cvd{\elacvd\global\havecvdtrue}
\let\siamproof=\proof

\def\proof{\siamproof\global\havecvdfalse}

\begin{document}



\vbox to 0pt{%
\vss%
\noindent%
Published in \textit{Electronic Journal of Linear Algebra}, \textbf{27}, 172--189, 2014.\\
\texttt{http://www.math.technion.ac.il/iic/ela/27.html}}

\title{On Kronecker Quotients}

\author{
 Yorick Hardy\thanks{
 Department of Mathematical Sciences,
 University of South Africa,
 Johannesburg, South Africa (hardyy@unisa.ac.za, yorickhardy@gmail.com).}}

\pagestyle{myheadings}
\markboth{Y.\ Hardy}{On Kronecker Quotients}
\maketitle

\begin{abstract}
 Leopardi introduced the notion of a Kronecker quotient in
 [Paul Leopardi. A generalized {FFT} for {C}lifford algebras.
 {\em Bulletin of the Belgian Mathematical Society}, 11:663--688, 2005.].
 This article considers the basic properties that a Kronecker quotient should
 satisfy and additional properties which may be satisfied.
 A class of Kronecker quotients for which
 these properties have a natural description is completely characterized.
 Two examples of types of Kronecker quotients are described.
\end{abstract}


\section{Introduction}
\label{sec:intro}

Let $\mathcal{M}(F,m,n)$ denote the vector space of $m\times n$ matrices
over the field $F$, where $m,n\in\mathbb{N}$.  Let
$\mathcal{M}_{\text{nz}}(F,m,n)=\mathcal{M}(F,m,n)\setminus\{0_{m\times n}\}$, where
$0_{m\times n}$ is the $m\times n$ zero matrix.

Let $A\in\mathcal{M}(F,m,n)$ and $B\in\mathcal{M}(F,s,t)$. The Kronecker product
\cite{graham81a,schacke04,kronecker2,vanloan00a}
$A\otimes B$ is the $ms\times nt$ matrix over $F$ with entries
\begin{equation}
 \label{eq:kronent}
 (A\otimes B)_{(i-1)s+p,(j-1)t+q} := (A)_{i,j}(B)_{p,q},
\end{equation}
where
 $i\in\{1,\ldots,m\}$,
 $j\in\{1,\ldots,n\}$,
 $p\in\{1,\ldots,s\}$ and
 $q\in\{1,\ldots,t\}$.
In other words, we can write in block matrix form
\begin{equation*}
A\otimes B=\begin{bmatrix}
                    (A)_{1,1}B & (A)_{1,2}B & \ldots & (A)_{1,n} B \\
                    (A)_{2,1}B & (A)_{2,2}B & \ldots & (A)_{2,n} B \\
                    \vdots     & \vdots     & \ddots & \vdots      \\
                    (A)_{m,1}B & (A)_{m,2}B & \ldots & (A)_{m,n} B
           \end{bmatrix}.
\end{equation*}

From the block matrix form it is obvious that $B$ can be determined
from $A$ and $A\otimes B$ provided there exists an entry $(A)_{i,j}\neq 0$
in $A$ or equivalently $\|A\|\neq 0$ for some norm $\|\cdot\|$ in the
vector space $\mathcal{M}(F,m,n)$.

Methods for determining $B$ from $A$ and $A\otimes B$ have been described in the literature.
In \cite{vanloan93a}, van Loan and Pitsianis describe how to find $C$
(for given $A$ and $M$) which minimizes the Frobenius norm $\|A\otimes C - M\|_F$.
When $M=A\otimes B$, their technique yields $B$. This method is
the example given in Section \ref{sec:frobnorm}. The software package
QuCalc \cite{qucalc} provides \verb|krondiv| which can calculate
$B$ from $A\otimes B$ provided $A$ and $B$ are column vectors. The software
package Scilab \cite{scilab} provided the right and left Kronecker division operators \verb|./.|
and \verb|.\.| (earlier documentation appears in \cite{scilabarchive}).
In \cite{leopardi05a} Leopardi defined a Kronecker quotient $\lkq$ and
proved that from his definition $A\lkq(A\otimes B)=B$. This article
initiates an exploration of all Kronecker quotients and their properties.

The rest of the article is arranged as follows.
Section \ref{sec:quotient} introduces definitions and properties relating to the
Kronecker quotient.
In Section \ref{sec:linear}, we consider linear Kronecker
quotients,
in particular this allows us to characterize a class
of Kronecker quotients, the uniform Kronecker quotients, in Section \ref{sec:uniform}.
In Section \ref{sec:eg}, we give examples of uniform Kronecker quotients
provided by weighted averages and partial norms.
Section \ref{sec:conc} concludes the article with some open problems.

In this article, we denote by
$\{\,\mathbf{e}_{1,m},\,\mathbf{e}_{2,m},\,\ldots,\,\mathbf{e}_{m,m}\,\}$
the standard basis in the vector space $F^m$ (column vectors),
by $I_m$ the $m\times m$ identity matrix, and by $A^T$ the transpose
of the matrix $A$.

\section{Kronecker quotients}
\label{sec:quotient}

Here we first consider the essential properties that a left
Kronecker quotient should satisfy, and then consider examples
which satisfy these properties such as the quotient defined
by Leopardi. Then we consider properties of Kronecker products
and the corresponding properties of Kronecker quotients.

\begin{ldefinition}{def:lkq}%
 Let 
 \begin{equation*}
  \lkq=\{\lkq_{m,n,s,t}:\mathcal{M}_{\text{nz}}(F,m,n)\times\mathcal{M}(F,ms,nt)
               \to\mathcal{M}(F,s,t),\;m,n,s,t\in\mathbb{N}\}
 \end{equation*}
 and define $A\lkq M:=\lkq_{m,n,s,t}(A,M)$
 for all $A\in\mathcal{M}_{\text{nz}}(F,m,n)$ and $M\in\mathcal{M}(F,ms,nt)$.
 If
 \begin{equation*}
  A\lkq(A\otimes B) = B
 \end{equation*}
 for all $m,n,s,t\in\mathbb{N}$, $A\in\mathcal{M}_{\text{nz}}(F,m,n)$ and $B\in\mathcal{M}(F,s,t)$,
 then $\lkq$ is a \emph{left Kronecker quotient}.
\end{ldefinition}

\begin{ldefinition}{def:rkq}%
 Let 
 \begin{equation*}
  \rkq=\{\rkq_{m,n,s,t}:\mathcal{M}(F,ms,nt)\times\mathcal{M}_{\text{nz}}(F,s,t)
               \to\mathcal{M}(F,m,n),\;m,n,s,t\in\mathbb{N}\}
 \end{equation*}
 and $M\rkq B:=\rkq_{m,n,s,t}(M,B)$
 for all $B\in\mathcal{M}_{\text{nz}}(F,s,t)$ and $M\in\mathcal{M}(F,ms,nt)$. If
 \begin{equation*}
  (A\otimes B)\rkq B = A
 \end{equation*}
 for all $m,n,s,t\in\mathbb{N}$, $A\in\mathcal{M}(F,m,n)$ and $B\in\mathcal{M}_{\text{nz}}(F,s,t)$,
 then $\rkq $ is a \emph{right Kronecker quotient}.
\end{ldefinition}

In \cite{vanloan93a}, van Loan and Pitsianis describe how to find $B$
(for given $A$ and $M$) which minimizes the Frobenius norm $\|A\otimes B - M\|_F$:
\begin{equation*}
 (B)_{i,j}=\frac{\tr(\tilde M_{ij}^TA)}{\|A\|_F^2}
\end{equation*}
where
\begin{equation*}
 \tilde
 M_{ij}:=\left(I_m\otimes\mathbf{e}_{i,s}\right)^T M
         \left(I_n\otimes\mathbf{e}_{j,t}\right).
\end{equation*}
This method provides a Kronecker quotient, which is the example given in Section \ref{sec:frobnorm}.
Leopardi defined a left Kronecker quotient in \cite{leopardi05a}:
\begin{equation*}
 A\lkq_{L}M
  :=\frac{1}{\nnz(A)}\sum_{(i,j)\in\nz(A)}
    \frac{M_{i,j}}{(A)_{i,j}}
\end{equation*}
where
\begin{equation*}
 M_{i,j}:=\left(\mathbf{e}_{i,m}\otimes I_s\right)^T M
          \left(\mathbf{e}_{j,n}\otimes I_t\right)
\end{equation*}
is the $s\times t$ matrix in the $i$-th row and $j$-th column
of the block structured matrix $M$ over $\mathcal{M}(F,m,n)\otimes\mathcal{M}(F,s,t)$ and
\begin{equation*}
 \nz(A):=
  \Big\{(i,j)\in\{1,2,\ldots m\}\times\{1,2,\ldots,n\}\,:\,(A)_{i,j}\neq0\Big\}
\end{equation*}
and  ${\nnz}(A):=|{\nz}(A)|$ is the number of non-zero entries in $A$.
He showed that $\lkq_L$ satisfies Definition \ref{def:lkq}.
Leopardi also showed that for $F=\reals$,
and $A\in\mathcal{M}_{\text{nz}}(\reals,2^n,2^n)$,
$\nnz(A)=2^n$ and all $C\in\mathcal{M}(\reals,2^n,2^n)$
\begin{equation*}
 A\lkq_L(C\otimes B)
  =(A'\bullet C)B,
\end{equation*}
where
\begin{equation*}
 (A')_{j,k}=
  \begin{cases}
   1/(A_{j,k})&(A)_{j,k}\neq 0\\
   0          &(A)_{j,k}=0
  \end{cases}
\end{equation*}
and $\bullet$ denotes the normalized Frobenius inner product.
This result is generalized in Section \ref{sec:uniform}.

\begin{ldefinition}{def:remainder}%
 Let $\lkq$
 denote a left Kronecker quotient, $A\in \mathcal{M}_{\text{nz}}(F,m,n)$ and
 $M\in\mathcal{M}(F,ms,nt)$.
 The $s\times t$ matrix $A\lkq M$ is the \emph{left Kronecker quotient
 of $M$ and $A$} and the $ms\times st$ matrix
 \begin{equation*}
  M \mathbin{\rem_L} A:=M-A\otimes(A\lkq M)
 \end{equation*}
 is the \emph{left Kronecker remainder of $M$ and $A$}
 with respect to $\lkq$.
 If $M \mathbin{\rem_L} A=0_{m\times n}$ then $A$ is a
 \emph{left Kronecker divisor of $M$}.
\end{ldefinition}

An analogous definition holds for the right Kronecker remainder.
In the remainder of the article we will only consider left Kronecker
quotients and remainders, analogous definitions and results for the
right Kronecker quotient and remainder are straightforward.

Many properties of the Kronecker product are described, for example, in
\cite{schacke04,kronecker2,vanloan00a}. These properties should be considered
when defining a Kronecker quotient. Some of the properties are
\begin{equation}
 \tag{K1}\label{eq:K1}
 (A\otimes B)^T=A^T\otimes B^T
\end{equation}
\begin{equation}
 \tag{K2}\label{eq:K2}
 \begin{array}{c}
   (A+C)\otimes B=A\otimes B+C\otimes B,\quad
   A\otimes(B+D)=A\otimes B+A\otimes D,\\
   k(A\otimes B)=(kA)\otimes B=A\otimes(kB)
   \end{array}
\end{equation}
\begin{equation}
 \tag{K3}\label{eq:K3}
 A\otimes (B\otimes C) = (A\otimes B)\otimes C
\end{equation}
\begin{equation}
 \tag{K4}\label{eq:K4}
 (A\otimes B)(G\otimes H)=(AG)\otimes(BH),
\end{equation}
where the matrices $C$ and $D$ are $m\times n$ and $s\times t$ respectively
and $G$ and $H$ are assumed to be compatible with $A$ and $B$ for the matrix
product, and $k\in F$. For $m=n$ and $s=t$ we have
\begin{equation}
 \tag{K5}\label{eq:K5}
 \tr(A\otimes B)=\tr(A)\;\tr(B)
\end{equation}
\begin{equation}
 \tag{K6}\label{eq:K6}
 \det(A\otimes B)=(\det A)^s\,(\det B)^m.
\end{equation}
We propose properties for Kronecker quotients corresponding to 
\eqref{eq:K1} -- \eqref{eq:K3} for Kronecker products.
Let $m,n,p,q,s,t,u,v\in\mathbb{N}$, $k\in F\setminus\{0\}$,
    $A\in\mathcal{M}_{\text{nz}}(F,m,n)$,
    $B\in\mathcal{M}_{\text{nz}}(F,s,t)$,
    $C\in\mathcal{M}_{\text{nz}}(F,n,p)$,
    $D\in\mathcal{M}_{\text{nz}}(F,t,q)$,
    $M,M_1,M_3\in\mathcal{M}(F,ms,nt)$,
    $M_2\in\mathcal{M}(F,nt,pu)$
    and $M_4\in\mathcal{M}(F,nt,qv)$.
\begin{equation}
 \tag{Q1}\label{eq:Q1}
 (A\lkq M)^T=A^T\lkq M^T\qquad (M\rkq B)^T=M^T\rkq B^T
\end{equation}
\begin{subequations}
 \begin{equation}
 \tag{Q2a}\label{eq:Q2a}
 \begin{array}{c}
  A\lkq (M_1+M_2)=A\lkq M_1 + A\lkq M_2,\quad
  A\lkq (kM) = k(A\lkq M),\\
  (M_1+M_2)\rkq B=M_1\rkq B+M_2\rkq B,\quad
  (kM)\rkq B = k(M\rkq B)
 \end{array}
 \end{equation}
 \begin{equation}
   \tag{Q2b}\label{eq:Q2b}
   (kA)\lkq M = \frac1k(A\lkq M),\quad
   M\rkq (kB) = \frac1k(M\rkq B)
 \end{equation}
\end{subequations}
\begin{equation}
\tag{Q3}\label{eq:Q3}
 A\lkq (B\lkq M) = (B\otimes A)\lkq M,\quad
 (M\rkq B)\rkq A = M\rkq(B\otimes A).
\end{equation}
For \eqref{eq:K4} -- \eqref{eq:K6} we propose the properties
\eqref{eq:Q4} -- \eqref{eq:Q6}.
\begin{equation}
\tag{Q4'}\label{eq:Q4}
(A\lkq M_1)(C\lkq M_2)=(AC)\lkq (M_1M_2)\atop
  (M_3\rkq B)(M_4\rkq D)=(M_3M_4)\rkq (BD)
\end{equation}
and for $m=n$ and $s=t$ we propose
\begin{equation}
 \tag{Q5'}\label{eq:Q5}
 \tr M=\tr(A)\;\tr(A\lkq M) =\tr(B)\;\tr(M\rkq B)
\end{equation}
\begin{equation}
 \tag{Q6'}\label{eq:Q6}
 \det(M)=(\det A)^s(\det(A\lkq M))^m =(\det B)^m(\det(M\rkq B))^s.
\end{equation}

Note that Definition \ref{def:lkq} does not imply \eqref{eq:Q1} -- \eqref{eq:Q6},
but \eqref{eq:Q1} -- \eqref{eq:Q6} do hold when the right hand arguments
of $\lkq$ are Kronecker products of an appropriate form, i.e.,
by Definition \ref{def:lkq},
\begin{equation}
 \tag*{Q1:}
 (A\lkq(A\otimes B))^T=B^T=(A^T\lkq (A\otimes B)^T)
\end{equation}
\begin{subequations}
 \begin{equation}
 \tag*{Q2a:}
 \begin{array}{c}
  A\lkq (A\otimes B_1+A\otimes B_2)
  =A\lkq (A\otimes(B_1+B_2))=B_1+B_2\\
  A\lkq(A\otimes B_1) + A\lkq(A\otimes B_2)=B_1+B_2\\
  A\lkq (k(A\otimes B)) = A\lkq(A\otimes kB)
   =kB = k(A\lkq(A\otimes B))
 \end{array}
 \end{equation}
 \begin{equation}
   \tag*{Q2b:}
   (kA)\lkq(A\otimes B) = (kA)\lkq\left(kA\otimes\frac1kB\right)
                        =\frac1k B
 \end{equation}
\end{subequations}
\begin{equation}
\tag*{Q3:}
 A\lkq (B\lkq (B\otimes A\otimes C))
  = A\lkq (A\otimes C) = C
  = (B\otimes A)\lkq (B\otimes A\otimes C)
\end{equation}
\begin{subequations}
\begin{equation}
 \tag*{Q4':}
 \begin{array}{c}
  (A_1\lkq (A_1\otimes B_1))(A_2\lkq (A_2\otimes B_2))
   = B_1B_2\\
  (A_1A_2)\lkq((A_1\otimes B_1)(A_2\otimes B_2))
   = (A_1A_2)\lkq((A_1A_2)\otimes(B_1B_2)) =B_1B_2
 \end{array}
\end{equation}
\end{subequations}
when $A_1A_2\neq 0$. For $m=n$ and $s=t$ 
\begin{equation}
 \tag*{Q5':}
 \tr (A\otimes B)=\tr A \tr B=\tr(A)\;\tr(A\lkq (A\otimes B))
\end{equation}
\begin{equation}
 \tag*{Q6':}
 \det(A\otimes B)=(\det A)^s(\det B)^m=(\det A)^s(\det(A\lkq (A\otimes B))^m.
\end{equation}
As a consequence, the map $M\mapsto A\lkq M$ is linear on
$\{A\otimes B\,|\,B\in\mathcal{M}(F,s,t)\}$
for all $s,t\in\mathbb{N}$.

The property \eqref{eq:Q4}
does not hold for any Kronecker quotient (as a counter
example choose non-zero $A,C$ such that $AC=0$),
we also find that \eqref{eq:Q5} cannot be satisfied for all $M$ and $A$
by any Kronecker quotient (choosing $M$ and $A$ with $\tr M\neq 0$ and $\tr A=0$
provides the counter example) and similarly for \eqref{eq:Q6}. Consequently,
we primarily consider \eqref{eq:Q1}, \eqref{eq:Q2a}, \eqref{eq:Q2b} and \eqref{eq:Q3} in this article.
A restriction of \eqref{eq:Q5} will also be considered.

\section{Linear Kronecker quotients}
\label{sec:linear}

For Kronecker quotients $\lkq$ satisfying \eqref{eq:Q2a} (i.e., linear
in the right argument), we may express the Kronecker quotient
in terms of a set of linear operators $q_{A,s,t}$,
\begin{equation*}
 A\in\mathcal{M}_{\textrm{nz}}(F,m,n),\,
   M\in\mathcal{M}(F,ms,nt)
   \quad\Rightarrow\quad A\lkq M=q_{A,s,t}(M),
\end{equation*}
where
\begin{equation*}
 \begin{split}
  \{\,q_{A,s,t}:\mathcal{M}(F,ms,nt)\to\mathcal{M}(F,s,t)\,\,:\,\,
  m,n,s,t\in\mathbb{N},\,
  A\in\mathcal{M}(F,m,n),\,\\
  \forall B\in\mathcal{M}(F,s,t)\,\,\,q_{A,s,t}(A\otimes B)=B\,\}
 \end{split}
\end{equation*}
is a set of full rank linear operators. Every
$M\in\mathcal{M}(F,ms,nt)$
can be written uniquely in the form $M=A\otimes B + R$, where
$B\in\mathcal{M}(F,s,t)$ and $R\in\ker(q_{A,s,t})$. Thus, we
have the following characterization of \eqref{eq:Q1} -- \eqref{eq:Q3}
for linear Kronecker quotients.

\begin{ltheorem}{thm:linprop}
 A linear Kronecker quotient $\lkq$
 \begin{itemize}
  \item[(i)]
   satisfies \eqref{eq:Q1} if and only if
   $R\in\ker(q_{A,s,t})\Rightarrow R^T\in\ker(q_{A^T,t,s})$
  \item[(ii)]
   satisfies \eqref{eq:Q2b} if and only if
   $R\in\ker(q_{A,s,t})\Rightarrow R\in\ker(q_{kA,s,t})$\\
   (or equivalently $q_{kA,s,t}(M)=\frac1k q_{A,s,t}(M)$
    for all $M\in\mathcal{M}(F,ms,nt)$)
  \item[(iii)]
   satisfies \eqref{eq:Q3} if and only if
   $q_{A,p,q}(q_{B,mp,nq}(M'))=q_{B\otimes A,p,q}(M')$
 \end{itemize}
 for all $m,n,p,q,s,t\in\mathbb{N}$, $k\in F\setminus\{0\}$,
 $A\in\mathcal{M}_{\textrm{nz}}(F,m,n)$,
 $M'\in\mathcal{M}(F,msp,ntq)$,
 and
 $M=A\otimes B+R\in\mathcal{M}(F,ms,nt)$ such that
 $R\in\ker(q_{A,s,t})$.
\end{ltheorem}
\begin{proof}
 \begin{itemize}
  \item[\textit{(i)}]
   Suppose $(A\lkq M)^T=A^T\lkq M^T$, then inserting $M=A\otimes B+R$
   into \eqref{eq:Q1} yields
   \begin{equation*}
    (A\lkq M)^T = [q_{A,s,t}(A\otimes B+R)]^T=B^T=A^T\lkq(A^T\otimes B^T+R^T),
   \end{equation*}
   since $q_{A,s,t}$ is linear and $R\in\ker(q_{A,s,t})$, and by linearity of $q_{A^T,t,s}$,
   \begin{equation*}
    A^T\lkq(A^T\otimes B^T+R^T)=q_{A^T,t,s}(A^T\otimes B^T+R^T)=B^T+q_{A^T,t,s}(R^T).
   \end{equation*}
   The above two equations show that $R^T\in\ker(q_{A^T,t,s})$.
   Conversely, suppose that $R\in\ker(q_{A,s,t})$
   implies $R^T\in\ker(q_{A^T,t,s})$ then from $R\in\ker(q_{A,s,t})$ and linearity
   of $q_{A^T,t,s}$
   \begin{equation*}
    A^T\lkq M^T = q_{A^T,t,s}(A^T\otimes B^T+R^T)
                = B^T = [q_{A,s,t}(A\otimes B+R)]^T=(A\lkq M)^T.
   \end{equation*}
  \item[\textit{(ii)}]
   Since, for $k\neq 0$,
   \begin{align*}
    \frac1k(A\lkq M)&=\frac1kq_{A,s,t}(A\otimes B+R)
                     =\frac1kB\\
    (kA)\lkq(M)&=q_{kA,s,t}(A\otimes B+R)
                =q_{kA,s,t}((kA)\otimes(\textstyle\frac1k)B)
                +q_{kA,s,t}(R)\\
               &=\frac1kB+q_{kA,s,t}(R),
   \end{align*}
   equation \eqref{eq:Q2b} holds if and only if
   $R\in\ker(q_{kA,s,t})$, or equivalently
   \begin{equation*}
    q_{kA,s,t}(M)=\frac1kq_{A,s,t}(M).
   \end{equation*}
  \item[\textit{(iii)}]
   This follows directly when expressing $\lkq$ in terms of
   $\{q_{A,s,t}\}$. \cvd
 \end{itemize}
\end{proof}

\subsection{Partial Frobenius product}
\label{sec:partfrob}

A straightforward extension of the Frobenius inner product is useful
in the discussion of linear Kronecker quotients.

\begin{ldefinition}{}
 Let $m,n,s,t\in\mathbb{N}$, $A\in\mathcal{M}(F,m,n)$
 and $M\in\mathcal{M}(F,ms,nt)$. The \emph{partial Frobenius product}
 \begin{equation*}
 \circ:(\mathcal{M}(F,m,n)\times\mathcal{M}(F,ms,nt))\cup
        (\mathcal{M}(F,ms,nt)\times\mathcal{M}(F,m,n))\to\mathcal{M}(s,t)
 \end{equation*}
 is given entry wise by
 \begin{equation*}
  (A\circ M)_{u,v} = \sum_{j=1}^m\sum_{k=1}^n (A)_{j,k} (M)_{(j-1)s+u,(k-1)t+v}
 \end{equation*}
 and $M\circ A = A\circ M.$
\end{ldefinition}

Notice that the partial Frobenius product is not an inner product,
and is not associative, but is bilinear. When $F=\mathbb{R}$ and
$s=t=1$ the above definition reduces to the usual Frobenius inner product.
The following theorem is a consequence of this definition.

\begin{ltheorem}{thm:frob}%
 \begin{itemize}
  \item[(i)]  If $A\circ B$ is defined then $(A\circ B)^T=A^T\circ B^T$.
  \item[(ii)] If $A\in\mathcal{M}(F,m,n)$, $B\in\mathcal{M}(F,ms,nt)$ and $C\in\mathcal{M}(F,p,q)$,\\
              then $A\circ(B\otimes C)=(A\circ B)\otimes C$.
  \item[(iii)] If $A\in\mathcal{M}(F,msp,ntq)$, $B\in\mathcal{M}(F,s,t)$ and $C\in\mathcal{M}(F,p,q)$,\\
               then $A\circ(B\otimes C)=(A\circ B)\circ C$.
  \item[(iv)] If $A\in\mathcal{M}(F,m,n)$ then $A\circ (\mathbf{e}_{j,m}\mathbf{e}_{k,n}^T)=(A)_{j,k}$.
  \item[(v)] If $A\in\mathcal{M}(F,n,n)$ then $A\circ I_n=\tr A$.
 \end{itemize}
\end{ltheorem}
\begin{proof}
 \textit{(i)}, \textit{(iv)} and \textit{(v)} follow straightforwardly from the definition. For \textit{(ii)},
 since $B\otimes C\in\mathcal{M}(F,msp,ntq)$
 and every $u\in\{1,\ldots,sp\}$ can be written in the form
 $u=(u_1-1)p+u_2$, where $u_1\in\{1,\ldots,s\}$ and $u_2\in\{1,\ldots,p\}$
 (and similarly for $v$) and using \eqref{eq:kronent}
 \begin{align*}
  \lefteqn{(A\circ(B\otimes C))_{u,v}}\qquad&\\
    &=\sum_{j=1}^m\sum_{k=1}^n (A)_{j,k} (B\otimes C)_{(j-1)sp+u,(k-1)tq+v}\\
    &=\sum_{j=1}^m\sum_{k=1}^n (A)_{j,k} (B\otimes C)_{(j-1)sp+(u_1-1)p+u_2,(k-1)tq+(v_1-1)q+v_2}\\
    &=\sum_{j=1}^m\sum_{k=1}^n (A)_{j,k} (B\otimes C)_{((j-1)s+u_1-1)p+u_2,((k-1)t+v_1-1)q+v_2}\\
    &=\sum_{j=1}^m\sum_{k=1}^n (A)_{j,k} (B)_{(j-1)s+u_1,(k-1)t+v_1}(C)_{u_2,v_2}\\
    &=(A\circ B)_{u_1,v_1}(C)_{u_2,v_2}=((A\circ B)\otimes C)_{u,v}.
 \end{align*}
 For \textit{(iii)} we have (as above, using the fact that $j\in\{1,\ldots,sp\}$
 can be written in the form $j=(j_1-1)p+j_2$ and similarly for $k$)
 \begin{align*}
  \lefteqn{(A\circ(B\otimes C))_{u,v}}\qquad&\\
     &=\sum_{j=1}^{sp}\sum_{k=1}^{tq} (A)_{(j-1)m+u,(k-1)n+v}(B\otimes C)_{j,k}\\
     &=\sum_{j_1=1}^s\sum_{j_2=1}^p\sum_{k_1=1}^t\sum_{k_2=1}^q
       (A)_{((j_1-1)p+j_2-1)m+u,((k_1-1)q+k_2-1)n+v}(B)_{j_1,k_1}(C)_{j_2,k_2}\\
     &=\sum_{j_1=1}^s\sum_{j_2=1}^p\sum_{k_1=1}^t\sum_{k_2=1}^q
       (A)_{(j_1-1)pm+(j_2-1)m+u,(k_1-1)nq+(k_2-1)n+v}(B)_{j_1,k_1}(C)_{j_2,k_2}\\
     &=\sum_{j_2=1}^p\sum_{k_2=1}^q
       (A\circ B)_{(j_2-1)m+u,(k_2-1)n+v}(C)_{j_2,k_2}\\
     &=(A\circ B)\circ C. \cvd
 \end{align*}
\end{proof}

When $s=t=1$, \textit{(ii)} simplifies to
\begin{equation*}
 A\circ(B\otimes C)=(A\circ B) C.
\end{equation*}

\subsection{Partial Frobenius product and linear Kronecker quotients}
\label{sec:partfroblkq}

Let $M\in\mathcal{M}(F,ms,nt)$, then $M$ can be written in the form
\begin{equation*}
 M=\sum_{j=1}^m\sum_{k=1}^n\sum_{u=1}^s\sum_{v=1}^t
   (M)_{(j-1)s+u,(k-1)t+v} \mathbf{e}_{j,m}\mathbf{e}_{k,n}^T\otimes\mathbf{e}_{u,s}\mathbf{e}_{v,t}^T.
\end{equation*}
For linear Kronecker quotients $\lkq$
described by the set of linear maps $\{q_{A,s,t}\}$ we have
\begin{equation*}
 A\lkq M=\sum_{j=1}^m\sum_{k=1}^n\sum_{u=1}^s\sum_{v=1}^t
          (M)_{(j-1)s+u,(k-1)t+v} q_{A,s,t}(\mathbf{e}_{j,m}\mathbf{e}_{k,n}^T\otimes\mathbf{e}_{u,s}\mathbf{e}_{v,t}^T)
\end{equation*}
and defining the $(st)^2$ matrices $Q_{u,v,u',v'}(A)\in\mathcal{M}(F,m,n)$ by
\begin{equation*}
 q_{A,s,t}(\mathbf{e}_{j,m}\mathbf{e}_{k,n}^T\otimes\mathbf{e}_{u,s}\mathbf{e}_{v,t}^T)
  =\sum_{u'=1}^s\sum_{v'=1}^t (Q_{u,v,u',v'}(A))_{j,k} \mathbf{e}_{u',s}\mathbf{e}_{v',t}^T
\end{equation*}
we find
\begin{align*}
 A\lkq M&=\sum_{u,u'=1}^s\sum_{v,v'=1}^t\left(\sum_{j=1}^m\sum_{k=1}^n
           (M)_{(j-1)s+u,(k-1)t+v}(Q_{u,v,u',v'}(A))_{j,k}\right)
             \mathbf{e}_{u',s}\mathbf{e}_{v',t}^T\\
        &=\sum_{u,u'=1}^s\sum_{v,v'=1}^t(Q_{u,v,u',v'}(A)\circ M)_{u,v}
             \mathbf{e}_{u',s}\mathbf{e}_{v',t}^T
\end{align*}
which can also be written as
\begin{equation}
\label{eq:decomp}
  A\lkq M =\sum_{u,u'=1}^s\sum_{v,v'=1}^t
    \mathbf{e}_{u',s}\mathbf{e}_{u,s}^T(Q_{u,v,u',v'}(A)\circ M)\mathbf{e}_{v,t}\mathbf{e}_{v',t}^T.
\end{equation}
When $M=A\otimes B$ we must have (applying Theorem \ref{thm:frob} \textit{(ii)})
\begin{equation*}
 A\lkq(A\otimes B)
  =\sum_{u,u'=1}^s\sum_{v,v'=1}^t
    (Q_{u,v,u',v'}(A)\circ A)
    (\mathbf{e}_{u',s}\mathbf{e}_{u,s}^T)B(\mathbf{e}_{v,t}\mathbf{e}_{v',t}^T)=B
\end{equation*}
for all $B\in\mathcal{M}(F,s,t)$. In particular, considering $B=\mathbf{e}_{x,s}\mathbf{e}_{y,t}^T$
for all $x\in\{1,\ldots,s\}$ and $y\in\{1,\ldots,t\}$ we obtain
\begin{equation}
 \label{eq:realize}
 Q_{u,v,u',v'}(A)\circ A = \delta_{u,u'}\delta_{v,v'}.
\end{equation}
It is straightforward to verify that \eqref{eq:realize} also implies $A\lkq (A\otimes B)=B$.
Thus, we have the following theorem.

\begin{ltheorem}{}
 Let $\lkq$ be a linear Kronecker quotient. Then for all
 $m,n,s,t\in\mathbb{N}$ and $A\in\mathcal{M}_{nz}(F,m,n)$ there
 exist $(st)^2$ matrices $Q_{u,u',v,v'}(A)\in\mathcal{M}(F,m,n)$, where
 $u,u'\in\{1,\ldots,s\}$, $v,v'\in\{1,\ldots,t\}$ and
 \begin{equation*}
  Q_{u,v,u',v'}(A)\circ A = \delta_{u,u'}\delta_{v,v'},
 \end{equation*}
 such that for all $M\in\mathcal{M}(F,ms,nt)$
 \begin{equation*}
  A\lkq M=
   \sum_{u,u'=1}^s\sum_{v,v'=1}^t
   \mathbf{e}_{u',s}\mathbf{e}_{u,s}^T(Q_{u,v,u',v'}(A)\circ M)\mathbf{e}_{v,t}\mathbf{e}_{v',t}^T.
 \end{equation*}
\end{ltheorem}

This result appears inconvenient, but provides the link to uniform Kronecker quotients.

\section{Uniform Kronecker quotients}
\label{sec:uniform}

For the case $s=t=1$ equation \eqref{eq:realize} becomes
$Q(A)\circ A=1$, where $Q(A):=Q_{1,1,1,1}(A)$ and $A\lkq M=Q(A)\circ M$.
This property will be used as the defining property for uniform
Kronecker quotients.

\begin{ldefinition}{}
 A left Kronecker quotient $\lkq$ is \emph{uniform} if for all
 $m,n\in\mathbb{N}$ and for every $A\in\mathcal{M}_{nz}(F,m,n)$ there
 exists  $Q(A)\in\mathcal{M}_{nz}(F,m,n)$ such that
 $Q(A)\circ A=1$ and for all $s,t\in\mathbb{N}$, $M\in\mathcal{M}(F,ms,nt)$
 \begin{equation*}
  A\lkq M=Q(A)\circ M.
 \end{equation*}
\end{ldefinition}

In other words, we have chosen the linear Kronecker quotient given by
\begin{equation*}
 Q_{u,v,u',v'}(A) = \delta_{u,u'}\delta_{v,v'}Q(A),\qquad Q(A)\circ A=1.
\end{equation*}
Similar to Leopardi's method, if $A\in\mathcal{M}_{nz}(F,m,n)$ and
$C\in\mathcal{M}(F,m,n)$ then
\begin{equation*}
 A\lkq(C\otimes B)=(Q(A)\circ C) B.
\end{equation*}
Note also that $Q(A)$ is uniquely determined by $\lkq$:
\begin{equation*}
 A\lkq(\mathbf{e}_{j,m}\mathbf{e}_{k,n}^T)=Q(A)\circ(\mathbf{e}_{j,m}\mathbf{e}_{k,n}^T)
  =(Q(A))_{j,k},
\end{equation*}
where we used Theorem \ref{thm:frob} \textit{(iv)}.

\begin{ldefinition}{}
 The matrices
 \begin{equation*}
  \left\{\,Q(A)\,:\,m,n\in\mathbb{N},A\in\mathcal{M}_{\text{nz}}(F,m,n)\,\right\}
 \end{equation*}
 are a called a \emph{realization} of a uniform Kronecker quotient $\lkq$, or equivalently
 $\lkq$ is realized by the matrices $Q(A)$.
\end{ldefinition}

\begin{ltheorem}{}
 \label{thm:realize}%
 A uniform left Kronecker quotient $\lkq$, realized by the matrices $Q(A)$
 satisfies \eqref{eq:Q2a}, and 
 \begin{itemize}
  \item[(i)] satisfies \eqref{eq:Q1} if and only if $Q(A^T)=(Q(A))^T$ and
  \item[(ii)] satisfies \eqref{eq:Q2b} if and only if $Q(kA)=\frac1{k} Q(A)$ and
  \item[(iii)] satisfies \eqref{eq:Q3} if and only if $Q(B\otimes A)=Q(B)\otimes Q(A)$
 \end{itemize}
 for all $k\in F\setminus\{0\}$, $A\in\mathcal{M}_{\text{nz}}(F,m,n)$
 and $B\in\mathcal{M}_{\text{nz}}(F,s,t)$.
\end{ltheorem}
\begin{proof}
 \begin{itemize}
  \item[\textit{(i)}] Using Theorem \ref{thm:frob} \textit{(iv)},
             we find that \textit{(i)} follows directly from the fact that $\lkq$ is linear and
             \begin{equation*}
              (A\lkq (\mathbf{e}_{j,m}\mathbf{e}_{k,n}^T))^T
               =A^T\lkq(\mathbf{e}_{j,m}\mathbf{e}_{k,n}^T)^T
             \end{equation*}
             if and only if
             \begin{equation*}
              (Q(A)\circ (\mathbf{e}_{j,m}\mathbf{e}_{k,n}^T))^T
               = Q(A^T)\circ(\mathbf{e}_{k,n}\mathbf{e}_{j,m}^T)
             \end{equation*}
             if and only if $(Q(A))_{j,k}=(Q(A^T))_{k,j}$.
             Conversely,
             \begin{equation*}
              A^T\lkq M^T=Q(A^T)\circ M^T=Q(A)^T\circ M^T=(Q(A)\circ M)^T=(A\lkq M)^T
             \end{equation*}
             by Theorem \ref{thm:frob} \textit{(i)}.
  \item[\textit{(ii)}]
             Here we use that $\lkq$ is linear and
             \begin{equation*}
              (kA)\lkq (\mathbf{e}_{j,m}\mathbf{e}_{k,n}^T)
               =\frac1k(A\lkq(\mathbf{e}_{j,m}\mathbf{e}_{k,n}^T))
             \end{equation*}
             if and only if
             \begin{equation*}
              Q(kA)\circ (\mathbf{e}_{j,m}\mathbf{e}_{k,n}^T)
               = \frac1k(Q(A)\circ(\mathbf{e}_{j,m}\mathbf{e}_{k,n}^T))
             \end{equation*}
             if and only if $(Q(kA))_{j,k}=\frac1k(Q(A))_{j,k}$.
             Conversely,
             \begin{equation*}
              (kA)\lkq M=Q(kA)\circ M=\left(\frac1kQ(A)\right)\circ M
                        =\frac1k(Q(A)\circ M)=\frac1k(A\lkq M),
             \end{equation*}
             since $\circ$ is bilinear.
  \item[\textit{(iii)}]
             Let $A\in\mathcal{M}(F,m,n)$, $B\in\mathcal{M}(F,s,t)$
             and $M\in\mathcal{M}(F,msp,ntq)$. We have
             \begin{equation*}
              (B\otimes A)\lkq M = A\lkq (B\lkq M)
             \end{equation*}
             if and only if
             \begin{equation*}
              Q(B\otimes A)\circ M = Q(A)\circ(Q(B)\circ M)
                                   = (M\circ Q(B))\circ Q(A)
                                   = M\circ(Q(B)\otimes Q(A))
             \end{equation*}
             by commutativity of $\circ$, Theorem \ref{thm:frob} \textit{(iii)} and \textit{(iv)},
             and considering the entries of $Q(B\otimes A)$ and $Q(B)\otimes Q(A)$ given by
             $M=(\mathbf{e}_{u,s}\mathbf{e}_{v,t}^T)\otimes(\mathbf{e}_{j,m}\mathbf{e}_{k,n}^T)$. \cvd
 \end{itemize}
\end{proof}

Corollary 7.6 in \cite{leopardi05a} is generalized as follows.

\begin{ltheorem}{}
 \label{thm:proj}
 If $\{A_1,\ldots,A_{mn}\}$ is a basis for $\mathcal{M}(F,m,n)$
 and $\lkq$ is a uniform left Kronecker quotient realized by
 the matrices $Q(A)$, then it holds that
 \begin{equation*}
  \text{for all $s,t\in\mathbb{N}$ and $M\in\mathcal{M}(F,ms,nt)$:~~}
  M=\sum_{j=1}^{mn} A_j\otimes(A_j\lkq M)
 \end{equation*}
 if and only if $Q(A_j)\circ A_k=\delta_{j,k}$.
\end{ltheorem}
\begin{proof}
 Since $\{A_1,\ldots,A_{mn}\}$ is a basis for $\mathcal{M}(F,m,n)$,
 $M$ can be written in the form
 \begin{equation}
  \label{eq:expand}
  M=\sum_{j=1}^{mn}A_j\otimes B_j,
 \end{equation}
 where $B_j\in\mathcal{M}(F,s,t)$ for $j\in\{1,\ldots,mn\}$.
 Thus,
 \begin{equation*}
  A_j\lkq M=\sum_{k=1}^{mn}Q(A_j)\circ (A_k\otimes B_k)
           =\sum_{k=1}^{mn}(Q(A_j)\circ A_k)B_k.
 \end{equation*}
 Consequently, if $Q(A_j)\circ A_k=\delta_{j,k}$ then
 $A_j\lkq M=B_j$ and
 \begin{equation*}
  M=\sum_{j=1}^{mn} A_j\otimes B_j = \sum_{j=1}^{mn}A_j\otimes(A_j\lkq M).
 \end{equation*}
 Conversely, suppose that for all $s,t\in\mathbb{N}$ and $M\in\mathcal{M}(F,ms,nt)$
 \begin{equation}
  \label{eq:toprove}
  M=\sum_{j=1}^{mn}A_j\otimes(A_j\lkq M).
 \end{equation}
 In particular, choose $st\geq mn$ and $M$ such that $\{B_1,\ldots,B_{mn}\}$ is linearly
 independent in $\mathcal{M}(F,s,t)$ so that (from \eqref{eq:expand}, and inserting
 \eqref{eq:expand} into \eqref{eq:toprove})
 \begin{equation*}
  M=\sum_{j=1}^{mn}A_j\otimes B_j=\sum_{j,k=1}^{mn}A_j\otimes (Q(A_j)\circ A_k)B_k,
 \end{equation*}
 which, since $\{A_1,\ldots,A_{mn}\}$ is a basis, yields
 \begin{equation*}
  B_j=\sum_{k=1}^{mn}(Q(A_j)\circ A_k)B_k,
 \end{equation*}
 and since $\{B_1,\ldots,B_{mn}\}$ was chosen to be linearly independent,
 \begin{equation*}
  Q(A_j)\circ A_k=\delta_{j,k}. \cvd
 \end{equation*}
\end{proof}

\subsection{Uniform Kronecker quotients and \eqref{eq:Q5}}
Suppose $A\in\mathcal{M}_{\text{nz}}(F,m,m)$ with $\tr(A)\neq 0$, $\lkq$
is a uniform Kronecker quotient and
\begin{equation*}
 \tr(M)=\tr(A)\;\tr(A\lkq M),
\end{equation*}
where $M=C\otimes B$ for some $C\in\mathcal{M}(F,m,m)$
and $B\in\mathcal{M}(F,s,s)$. It follows that $\tr(B)=0$ or
\begin{equation*}
 \tr(C) = \tr(A)(Q(A)\circ C)
\end{equation*}
which can be rewritten as (Theorem \ref{thm:frob} \textit{(v)})
\begin{equation*}
 I_{m}\circ C = \left[\tr(A)Q(A)\right]\circ C
\end{equation*}
and, by considering all $C$ from $\mathcal{M}(F,m,m)$
and $B\in\mathcal{M}(F,s,s)$,
\begin{equation*}
 Q(A)=\frac{1}{\tr(A)}I_m.
\end{equation*}

Thus, we have the following theorem.

\begin{ltheorem}{}
\label{thm:realizetrace}%
A uniform left Kronecker quotient $\lkq$, realized by the matrices $Q(A)$
satisfies \eqref{eq:Q5} restricted to $\tr(A)\neq 0$ if and only if
\begin{equation}
 \tag{TR}\label{eq:TR}
 Q(A)=\frac{1}{\tr(A)}I_m
\end{equation}
for all $A\in\mathcal{M}_{\text{nz}}(F,m,m)$ with $\tr(A)\neq 0$.
\end{ltheorem}

Notice that \eqref{eq:Q5} is incompatible with \eqref{eq:Q3},
for example, since
\begin{equation*}
 \begin{bmatrix}1&1&1&1\cr1&1&1&1\cr1&1&1&1\cr1&1&1&1\end{bmatrix}
 =\begin{bmatrix}1&1\cr1&1\end{bmatrix}\otimes\begin{bmatrix}1&1\cr1&1\end{bmatrix}
 =\begin{bmatrix}1\cr1\cr1\cr1\end{bmatrix}\otimes\begin{bmatrix}1&1&1&1\end{bmatrix}
\end{equation*}
we have by \eqref{eq:Q3} and \eqref{eq:TR}
\begin{align*}
 Q\left(\begin{bmatrix}1&1&1&1\cr1&1&1&1\cr1&1&1&1\cr1&1&1&1\end{bmatrix}\right)
  &=Q\left(\begin{bmatrix}1&1\cr1&1\end{bmatrix}\right)\otimes
    Q\left(\begin{bmatrix}1&1\cr1&1\end{bmatrix}\right)
  =\frac14 I_4 \\
\intertext{and by \eqref{eq:Q3}}
 Q\left(\begin{bmatrix}1&1&1&1\cr1&1&1&1\cr1&1&1&1\cr1&1&1&1\end{bmatrix}\right)
  &=Q\left(\begin{bmatrix}1\cr1\cr1\cr1\end{bmatrix}\right)\otimes
    Q\left(\begin{bmatrix}1&1&1&1\end{bmatrix}\right)
\end{align*}
which is impossible since $I_4$ has rank 4, while the last expression
has rank 1.

\section{Examples of uniform Kronecker quotients}
\label{sec:eg}

\subsection{Weighted average uniform Kronecker quotients}
\label{sec:wkq}

A generalization of Leopardi's method is as follows.

\begin{ldefinition}{}
Let $A\in\mathcal{M}_{\text{nz}}(F,m,n)$ and $M\in\mathcal{M}(F,ms,nt)$.
Let
\begin{equation*}
 W:=\{W_{m,n}:\mathcal{M}_{\text{nz}}(F,m,n)\to\mathcal{M}_{\text{nz}}(F,m,n),\,
       m,n\in\mathbb{N}\}
\end{equation*}
such that for all $A\in \mathcal{M}_{\text{nz}}(F,m,n)$
\begin{equation*}
 \sum_{(i,j)\in{\nz}(A)} (W_{m,n}(A))_{i,j} = 1.
\end{equation*}
The \emph{left weighted average Kronecker quotient} for the weights
$W$ is defined as
\begin{equation*}
A\lkq_W M:=\sum_{(i,j)\in{\nz}(A)}
     (W_{m,n}(A))_{i,j} \frac{M_{i,j}}{(A)_{i,j}},
\end{equation*}
where $M_{i,j}$ is the $s\times t$ matrix in the $i$-th row and $j$-th column
of the block structured matrix $M$ over $\mathcal{M}(F,m,n)\otimes\mathcal{M}(F,s,t)$.
\end{ldefinition}

The weighted average Kronecker quotient is uniform and is realized by
\begin{equation*}
 (Q(A))_{i,j}=
  \begin{cases}
   (W_{m,n}(A))_{i,j}/(A)_{i,j}&(A)_{i,j}\neq 0,\\
   0          &(A)_{i,j}=0.
  \end{cases}
\end{equation*}

Examples include (for a field $F$ with characteristic 0)
\begin{equation*}
 (W_{L}(A))_{i,j}=\frac{1-\delta_{(A)_{i,j},0}}{{\nnz}(A)}
\end{equation*}
for Leopardi's method and (for $F=\complex$ or $F=\reals$)
\begin{equation*}
 (W_F(A))_{i,j}=\frac{|(A)_{i,j}|^2}{\|A\|_F^2}
\end{equation*}
where $\|\cdot\|_F$ denotes the Frobenius norm.

\begin{ltheorem}{}
The properties \eqref{eq:Q2a} and \eqref{eq:Q2b} hold for $\lkq_{W}$.
The property \eqref{eq:Q1} holds for $\lkq_{W}$ if and only if
\begin{equation*}
 W_{n,m}(A^T)=\big(W_{m,n}(A)\big)^T.
\end{equation*}
The property \eqref{eq:Q3} holds for $\lkq_{W}$ if and only if
\begin{equation*}
 W_{ms,nt}(A\otimes B)=W_{m,n}(A)\otimes W_{s,t}(B)
\end{equation*}
for all $m,n,s,t\in\mathbb{N}$, $A\in\mathcal{M}_{\text{nz}}(F,m,n)$
and $B\in\mathcal{M}_{\text{nz}}(F,s,t)$.
\end{ltheorem}

The properties \eqref{eq:Q1}, \eqref{eq:Q2a}, \eqref{eq:Q2b} 
and \eqref{eq:Q3} hold for Leopardi's method and also for the
weighted average quotient defined by $W_F$.

\subsection{Partial norms}
Inspired by the partial trace (see for example \cite{carlen})
we define the notion of a partial norm. In this case, the underlying field is
$F=\complex$ or $F=\reals$. 

\begin{ldefinition}{}
Let $m,n,s,t\in\mathbb{N}$, $A\in\mathcal{M}(F,m,n)$, $\| \cdot \|$ denote a
norm on $\mathcal{M}(F,m,n)$ and let
$\|\cdot\|_{A,L}:\mathcal{M}(F,ms,nt)\to\mathcal{M}(F,s,t)$.
In other words the definition of $\|\cdot\|_{A,L}$ depends
on $A$ and $\|\cdot\|$. If for all $B\in\mathcal{M}(F,s,t)$,
\begin{equation*}
 \| A\otimes B\|_{A,L}:=\|A\| B
\end{equation*}
then $\|\cdot\|_{A,L}$ is a \emph{left partial norm}
with respect to $A$ and $\|\cdot\|$.
\end{ldefinition}

When $\|A\|\neq 0$ we may write
\begin{equation*}
 A\lkq (A\otimes B) = B = \frac{\| A\otimes B\|_{A,L}}{\|A\|}
\end{equation*}
which is sufficient to define $A\lkq$, i.e.,
\begin{equation*}
 A\lkq M:=\frac{\|M\|_{A,L}}{\|A\|}.
\end{equation*}

\begin{ltheorem}{}
 A left Kronecker quotient $\lkq$ defined by a left partial norm,
 given by the norm $\|\cdot\|$, is uniform (realized by the matrices $Q(A)$)
 if
 \begin{equation*}
  A\lkq M:=\frac{\|M\|_{A,L}}{\|A\|}
         =Q(A)\circ M,\qquad
  (Q(A))_{j,k}:=\frac{\|E_{j,k}\|_{A,L}}{\|A\|}
 \end{equation*}
 for all $m,n,s,t\in\mathbb{N}$, $A\in\mathcal{M}_{\text{nz}}(\complex,m,n)$ and
 $M\in\mathcal{M}(\complex,ms,nt)$,
 where $E_{j,k}\equiv E_{j,k}\otimes 1$ is the $m\times n$ matrix with the entry 1
 in the $j$-th row and $k$-th column and 0 for all other entries.
\end{ltheorem}

Next we discuss two examples, namely the uniform Kronecker
quotients given by the Frobenius norm and the operator norm.

\subsubsection{Example: Frobenius norm}
\label{sec:frobnorm}%

Let $A\in\mathcal{M}(\complex,m,n)$ be non-zero.
Consider the Frobenius norm (also known as the Hilbert-Schmidt norm)
\begin{equation*}
 \|A\|_F:=\sqrt{\tr(A^*A)},
\end{equation*}
where $A^*$ is the transposed complex conjugate of $A$.
Let $B\in\mathcal{M}(\complex,s,t)$.
Then
\begin{equation*}
  \frac1{\|A\|_F}\sum_{d=1}^n\left[A\mathbf{e}_{d,n}\otimes I_s\right]^*
  (A\otimes B)
  \left[\mathbf{e}_{d,n}\otimes I_t\right]
  =\|A\|_F B.
\end{equation*}

\begin{ldefinition}{}
The \emph{left Frobenius partial norm} of $M\in\mathcal{M}(\complex, ms,nt)$
with respect to $A\in\mathcal{M}(\complex,m,n)$ is defined as
\begin{equation*}
  \|M\|_{A,L,F}:=
  \frac1{\|A\|_F}\sum_{d=1}^n\left[A\mathbf{e}_{d,n}\otimes I_s\right]^*
  M
  \left[\mathbf{e}_{d,n}\otimes I_t\right].
\end{equation*}
\end{ldefinition}

\begin{ldefinition}{}
Let $M\in\mathcal{M}(\complex,ms,nt)$ and $A\in\mathcal{M}_{\text{nz}}(\complex,m,n)$.
The \emph{left Frobenius Kronecker quotient} is defined as
\begin{equation*}
 A\lkq_{F} M:=\frac{\|M\|_{A,L,F}}{\|A\|_F}.
\end{equation*}
\end{ldefinition}

The left Frobenius Kronecker quotient is uniform and is
realized by 
\begin{equation*}
 Q(A)=\frac{\overline{A}}{\|A\|_F^2}
\end{equation*}
where $\overline{A}$ is the complex conjugate of $A$.
This is identical to the second example of weighted average
Kronecker quotients in Section \ref{sec:wkq}.

The Kronecker quotient induced by the Frobenius partial
norm is equivalent to finding the nearest Kronecker product \cite{vanloan93a},
where $A\lkq_{F}M$ minimizes 
\begin{equation*}
 \|M-A\otimes(A\lkq_{F}M)\|_F.
\end{equation*}

\subsubsection{Example: Operator norm}

Let $A\in\mathcal{M}(\complex,m,n)$.
Consider the operator norm
\begin{equation*}
\|A\|_O:=\max \sigma(A),
\end{equation*}
where $\sigma(A)$ is the set of singular values of $A$.
Let $B\in\mathcal{M}(\complex,s,t)$.
Suppose $A=U\Sigma V^*$ is a singular value decomposition
of $A$, where $U$ is an $m\times m$ unitary matrix, $V$ is an
$n\times n$ unitary matrix,
\begin{equation*}
\left(\Sigma\right)_{u,v}
      =\left\{\begin{matrix}
                \delta_{u,v}\sigma_u&u,v\leq\min\{m,n\}\\
                0                   &\text{otherwise}
              \end{matrix}\right.
\end{equation*}
and $\sigma_1\geq\sigma_2\geq\cdots\geq\sigma_{\min\{m,n\}}\geq0$
are the singular values of $A$. Then $\|A\|_O=\sigma_1$
and
\begin{equation*}
  \left[\left(U\mathbf{e}_{1,m}\right)^*\otimes I_s\right]
  (A\otimes B)
  \left[\left(V\mathbf{e}_{1,n}\right)\otimes I_t\right]
  =\|A\|_OB.
\end{equation*}

Since the singular value decomposition of $A$ is not unique
in general, the algorithm for calculating the singular value
decomposition influences the properties of the partial norm
and Kronecker quotient.  We assume that the singular value
decomposition is determined uniquely (by an appropriate algorithm
for example) in the following.

\begin{ldefinition}{}
The \emph{left operator partial norm} of $M\in\mathcal{M}(\complex,ms,nt)$
with respect to $A\in\mathcal{M}(\complex,m,n)$ is defined as
\begin{equation*}
  \|M\|_{A,L,O}:=
  \left[\left(U\mathbf{e}_{1,m}\right)^*\otimes I_s\right]
  M
  \left[\left(V\mathbf{e}_{1,n}\right)\otimes I_t\right],
\end{equation*}
where $A=U\Sigma V^*$ is the singular value decomposition of $A$.
\end{ldefinition}

\begin{ldefinition}{}
Let $M\in\mathcal{M}(\complex,ms,nt)$ and $A\in\mathcal{M}_{\text{nz}}(\complex,m,n)$.
The \emph{left operator Kronecker quotient} is defined as
\begin{equation*}
 A\lkq_{O} M:=\frac{\|M\|_{A,L,O}}{\|A\|_O}.
\end{equation*}
\end{ldefinition}

The left operator Kronecker quotient is uniform and realized by the matrices $Q(A)$ 
(for the singular value decomposition $A=U\Sigma V^*$) given by
\begin{equation*}
 Q(A)=\frac{\overline{U}\mathbf{e}_{1,m}(V\mathbf{e}_{1,n})^T}{\|A\|_O}.
\end{equation*}

\section{Conclusion}
\label{sec:conc}

We have presented the basic properties of Kronecker quotients.
We completely characterized the uniform Kronecker quotients for which
these properties have a natural description. Two examples of types
of uniform Kronecker quotients were described.

Many interesting open problems remain, including:
\begin{enumerate}
 \item Is it possible to express \eqref{eq:decomp} in terms of uniform
       Kronecker quotients, i.e., can all Kronecker quotients be decomposed
       in terms of uniform Kronecker quotients?
 \item Is there a restricted form of \eqref{eq:Q3} such that
       (restricted) uniform Kronecker quotients may satisfy \eqref{eq:Q5}?
 \item To characterize multiplicative Kronecker quotients, i.e., to abandon
       \eqref{eq:Q2a} (and possibly \eqref{eq:Q2b}) in favor of \eqref{eq:Q4} and
       possibly \eqref{eq:Q6} (for example, restricted to the non-singular matrices).
\end{enumerate}

\section*{Acknowledgment}

The author acknowledges the anonymous referee for many improvements
in the presentation of this work, in particular the presentation of
Theorem \ref{thm:linprop}, and the
notation and terminology of the partial Frobenius product. The anonymous
referee also pointed out that a theorem very similar to Theorem \ref{thm:proj}
may hold, which was subsequently reformulated as Theorem \ref{thm:proj}.

The author is supported by the National Research Foundation (NRF),
South Africa. This work is based upon research supported by the National
Research Foundation. Any opinion, findings and conclusions or recommendations
expressed in this material are those of the author(s) and therefore the
NRF do not accept any liability in regard thereto.

\end{document}